\documentclass[10pt]{article}
\usepackage[cp1251]{inputenc}
\usepackage[english,russian]{babel}     
\usepackage{amsfonts}
\usepackage{amsthm} 
\usepackage{amssymb}
\usepackage{latexsym} 

\newtheorem{teo}{\quad Theorem}
\newtheorem{lem}{\quad Lemma}
\newtheorem{cor}{\quad Corollary}

\renewcommand{\refname}{Literature}

\begin{document}
\Large

\noindent UDC 517.5

\bigskip\noindent
{\bf S.О.~Chaichenko}
(Donbass State Pedagogical University, Slavyansk, Ukraine)

\bigskip\noindent
{\bf  CONVERGENCE OF  FOURIER SERIES ON THE SYSTEM \\ OF RATIONAL FUNCTIONS ON THE REAL  AXIS}

\bigskip
{\noindent We consider the systems of rational functions $\{\Phi_n(z)\}, ~n \in \mathbb{Z}$, defined by fixed set points ${\bf a}:=\{a_k\}_{k=0}^{\infty}, ~ (\mathop{\rm Im} a_k>0)$, ${\bf b}:=\{b_k\}_{k=1}^{\infty}, ~ (\mathop{\rm Im} b_k<0)$ and is orthonormal on the real axis $\mathbb{R}.$ We have obtained the compact form of analogue of Dirichlet kernels of these systems on the real axis $\mathbb{R}.$  Using obtained representation we investigate the problems of convergence in the spaces $L_p(\mathbb{R}),~ p> 1,$ and pointwise convergence of Fourier series on the systems $\{\Phi_n(t)\},~ n \in \mathbb{Z},$ provided that the sequences of poles of these systems satisfies certain restrictions. We have proved statements that are analogues of the classical Theorems of Jordan-Dirichlet and Dini-Lipschitz of convergence of  Fourier series on the trigonometric system.}
\bigskip

\bigskip

{ \bf AMS:} 46E30, 42A10,41A17,41A20,41A25,41A27, 41A30.

\bigskip

{\bf 1. Orthogonal system of rational functions on the real axis.}
Let  ${\bf a}:=\{a_k\}_{k=0}^{\infty}, ~ (\mathop{\rm Im} a_k>0)$ be arbitrary sequence of
complex numbers from the upper half plane $\mathbb C_+:=\{z\in\mathbb C : \mathop{\rm Im}z>0\}$.
Then
\begin{equation}\label{Phi-system}
    \Phi_0^+(z):=\frac{\sqrt{\mathop{\rm Im} a_0}}{z-\overline a_0},~
    \Phi_n^+(z):=\frac{\sqrt{\mathop{\rm Im} a_n}}{z-\overline a_n}B_n^+(z),~n=1,2,\ldots,
\end{equation}
where
$$
    B_0^+(z):=1,~B_n^+(z):=\prod_{k=0}^{n-1}\chi_k^+
    \frac{z-a_k}{z-\overline a_k},~
    \chi_k^+:=\frac{|1+a^2_k|}{1+a^2_k},~n=1,2,\ldots.
$$
--- $n$-Blaschke product with zeros at the points $a_k,$ $k=0,1,...n-1$.

The system of functions $\{\Phi^+_n(z)\}_0^\infty$ introduced by M.M. Dzhrbashyan \cite{Djrbashan_1974} similarly, as it did S.~Takenaka and F.~Malmquist  \cite{Takenaka, Malmquist} in the case of space Hardy $H_2$ in the unit circle. In particular, in the article \cite{Djrbashan_1974} was shown that the system $\{\Phi^+_n(z)\}_0^\infty$ is orthonormal on the real axis $\mathbb R,$ i.e.
$$
    \frac{1}{\pi} \int\limits_{-\infty}^\infty \Phi_n^+ (x)\overline{\Phi_m^+(x)}dx=
    \left\{\matrix{0, \hfill& n\not=m,\cr
                    1,\hfill& n=m, }\right.\quad n,m=0,1,2,\ldots.
$$

Let, further,  ${\bf b}:=\{b_k\}_{k=1}^{\infty}, ~ (\mathop{\rm Im} b_k<0)$ be arbitrary sequence of complex numbers from the lower half plane $\mathbb C_-:=\{z\in\mathbb C : \mathop{\rm Im}z<0\}$. Then
$$
    \Phi_1^-(z):=\frac{\sqrt{-\mathop{\rm Im} b_1}}{z-\overline b_1},~
    \Phi_n^-(z):=\frac{\sqrt{-\mathop{\rm Im} b_n}}{z-\overline b_n}B_n^-(z),~n=2,3,\ldots,
$$
where
$$
    B_1^-(z):=1,~B_n^-(z):=\prod_{k=1}^{n-1}\chi_k^-\frac{z-b_k}{z-\overline b_k},~
    \chi_k^- :=\frac{|1+b^2_k|}{1+b^2_k},~n=2,3,\ldots.
$$
--- $n$-Blaschke product with zeros at the points $b_k,$ $k=1,2,...n-1$.

\begin{lem} \label{L.1}
The system of functions $\{\Phi_n^-(z)\}_1^\infty$ is orthonormal system on the real axis $\mathbb{R},$ i.e.
$$
    \frac{1}{\pi} \int\limits_{-\infty}^\infty \Phi_n^- (x)\overline{\Phi_m^-(x)}dx=
    \left\{\matrix{0,\hfill& n\not=m,\cr
                    1,\hfill& n=m, }\right.\quad n,m=1,2,\ldots.
$$
\end{lem}

{\bf Proof.}  Note first that
$$
   \frac{1}{\pi} \int\limits_{-\infty}^\infty \Phi_n^- (x)\overline{\Phi_n^-(x)}dx= \frac{1}{\pi} \int\limits_{-\infty}^\infty |\Phi_n^- (x)|^2dx=
$$
$$
    =\frac{-\mathop{\rm Im} b_n}{\pi} \int\limits_{-\infty}^\infty \frac{dx}{|x-\bar{b}_n|^2}=1,
    \quad (n=1,2,\ldots).
$$

Assuming now that $1\le m \le n-1,~ (n \ge 2),$ consider the integrals
$$
    \frac{1}{\pi} \int\limits_{-\infty}^\infty \Phi_n^- (x) \overline{\Phi_m^- (x)}dx=
    \frac{\sqrt{\mathop{\rm Im} b_n \mathop{\rm Im} b_m}}{\pi} \int\limits_{-\infty}^\infty
    \frac{\prod\limits_{k=m+1}^{n-1}  (x-b_k) \chi_k^-}{\prod\limits_{k=m}^n (x-\bar{b}_k)}dx.
$$
The integrand in the right side of this equality is the rational fraction with poles at the points
$\overline {b}_k,~k=m,m+1,\ldots,n,$ from the upper half plane $\mathbb{C}_+$. It is clear that for $|z|\to \infty,$ $ z \in \mathbb{C},$
this fraction have the order ${\cal O}(|z|^{-2}),$ that ensures convergence of this integral. Therefore, using Cauchy integral theorem for the lower half plane $\mathbb{C}_-$, we get
$$
    \frac{1}{\pi} \int\limits_{-\infty}^\infty
    \Phi_{n}^-(x) \overline{\Phi_m^- (x)} dx=0.
$$

Finally, going to the conjugate values we sure that the last equality holds for arbitrary natural $n\not= m.$  Lemma is proved.

\begin{lem}\label{L.2} For arbitrary $z, \zeta\in\mathbb C,$ $z\not=\overline\zeta$ and $n=0,1,2,\ldots,$
$m=1,2,\ldots,$ the identities hold
\begin{equation}\label{Dzhrbashyan identity+}
    \sum_{k=0}^{n-1}\overline{\Phi_k^+(\zeta)}\Phi_k^+(z)=\frac{1}{2i(\overline\zeta-z)}\bigg[ 1-
    \overline{B_n^+(\zeta)}B_n^+(z)\bigg],
\end{equation}
\begin{equation}\label{Dzhrbashyan identity-}
    \sum_{k=1}^{m-1}\overline{\Phi_k^-(\zeta)}\Phi_k^-(z)=\frac{1}{2i(\overline\zeta-z)}
    \bigg[\overline{B_m^-(\zeta)}B_m^-(z)-1\bigg],
\end{equation}
 where $\sum_{k=0}^{-1}=\sum_{k=1}^{0}=0.$
\end{lem}

The identities (\ref{Dzhrbashyan identity+}) and (\ref{Dzhrbashyan identity-}) are analogues of well-known formula Christoffel-Darboux for orthogonal polynomials. These relations  is  known in the literature also as the  identities of M.M.~Dzhrbashyan. Proof of the identity (\ref{Dzhrbashyan identity+}) was obtained in the works \cite {Djrbashan_1974, Savchuk_CH_umg_2014}. In order to verify the equality (\ref{Dzhrbashyan identity-}) we will use the method proposed in \cite{Savchuk_CH_umg_2014}.

{\bf Proof.} Consider identity
\begin{equation}\label{identity}
    (1-x_1)+(1-x_2)x_1+(1-x_3)x_1x_2+\ldots+(1-x_{m-1})x_1x_2\cdots x_{m-2}=
$$
$$
    =1-x_1x_2\cdots x_{m-1},
\end{equation}
which is true for any complex numbers $x_1, x_2\ldots, x_{m-1}$, $m\in\mathbb N.$

Put
$$
x_k=\frac{(z-b_k)(\overline \zeta-\overline b_k)}{(z-\overline b_k)(\overline \zeta-b_k)}.
$$
Then
$$
    1-x_k=\frac{(z-\overline b_k)(\overline \zeta-b_k)-(z-b_k)
    (\overline \zeta-\overline b_k)}{(z-\overline b_k)(\overline \zeta-b_k)}=
    \frac{2i \mathop{\rm Im}b_k}{(z-\overline b_k)(\overline \zeta-b_k)}(\overline\zeta-z)
$$
and in accordance with (\ref{identity})
$$
    \sum_{k=1}^{m-1}(1-x_k)\prod_{j=1}^{k-1}x_j=
    \sum_{k=1}^{m-1}\frac{2i\mathop{\rm Im} b_k}{(z-\overline b_k)
    (\overline \zeta-b_k)}(\overline\zeta-z)
    \prod_{j=1}^{k-1}\frac{(z-b_j)(\overline \zeta-\overline b_j)}{(z-\overline b_j)
    (\overline \zeta-b_j)}=
$$
$$
    =1-\prod_{j=1}^{m-1}\frac{(z-b_j)(\overline \zeta-\overline b_j)}
    {(z-\overline b_j)(\overline \zeta-b_j)}.
$$

Dividing both sides of equality at $2i(\overline\zeta-z)$ and taking into account that
$$
    \prod_{j=1}^{k-1}\frac{(z-b_j)(\overline \zeta-\overline b_j)}
    {(z-\overline b_j)(\overline \zeta-b_j)}=
    \prod_{j=1}^{k-1}\chi_j^-\frac{z-b_j}{z-\overline b_j}
    \cdot\prod_{j=1}^{k-1}\overline{\chi}_j^{~-}\frac{\overline \zeta-\overline b_j}
    {\overline \zeta-b_j}=B_k^-(z)\overline{B_k^-(\zeta)}
$$
$b_k \in \mathbb{C}_-$ and
$$
    \sum_{k=1}^{m-1}\frac{\mathop{\rm Im} b_k}{(z-\overline b_k)(\overline \zeta-d_k)}
    \prod_{j=1}^{k-1}\frac{(z-b_j)(\overline \zeta-\overline b_j)}
    {(z-\overline b_j)(\overline \zeta-a_j)}=
$$
$$
    =-\sum_{k=1}^{m-1}\frac{\sqrt{-\mathop{\rm Im} b_k}}{z-\overline b_k} B_k(z)
    \frac{\sqrt{-\mathop{\rm Im} b_k}}{\overline\zeta-b_k} \overline{B_k(\zeta)}=
$$
$$
    =-\sum_{k=1}^{m-1}\Phi_k^-(z)\overline{\Phi_k^-(\zeta)},
$$
we obtain (\ref{Dzhrbashyan identity-}). Lemma is proved.

Let now
\begin{equation}\label{system}
    \Phi_n(z)=\cases{\Phi_n^+(z), & $n =0,1,2,\ldots,$ \cr
                     \Phi_{-n}^-(z), & $n =-1,-2,\ldots.$}
\end{equation}

\begin{lem}\label{L.3}
The system of functions $\{\Phi_n(z)\},~ n \in \mathbb{Z}$ is orthonormal system on the real axis $\mathbb{R},$ i.e.
\begin{equation}\label{ort_Omega}
    \frac{1}{\pi} \int\limits_{-\infty}^\infty \Phi_n (x) \overline{\Phi_m(x)}dx=
    \left\{\matrix{0,\hfill& n\not=m,\cr
                    1,\hfill& n=m, }\right.\quad n,m\in\mathbb Z.
\end{equation}
\end{lem}

{\bf Proof.}  If $n,m=0,1,2,\ldots,$ or $n,m=-1,-2,\ldots,$  then equality (\ref{ort_Omega}) follows from the facts of orthonormality of the systems of functions $\{\Phi_n^+(z)\}_0^\infty$ and $\{\Phi_k^-(z)\}_1^\infty$. Assume $n=-1,-2,\ldots,$ $m=0,1,2,\ldots$. Then
$$
    \frac{1}{\pi} \int\limits_{-\infty}^\infty \Phi_n (x) \overline{\Phi_m(x)}dx=
     \frac{1}{\pi} \int\limits_{-\infty}^\infty
    \Phi_{-n}^-(x) \overline{\Phi_m^+(x)} dx=
$$
$$
    =\frac{\sqrt{-\mathop{\rm Im} b_{-n} \mathop{\rm Im} a_m }}{\pi}\int\limits_{-\infty}^\infty
    \frac{1}{x-\overline b_{-n}}\prod_{j=1}^{-n-1}\chi_j^-\frac{x-b_j}{x-\overline b_j}\cdot
    \overline{\frac{1}{x-\overline a_{m}}\prod_{l=0}^{m-1}\chi_l^+ \frac{x-a_l}{x-\overline a_l}}~dx=
$$
\begin{equation}\label{0}
    =\frac{\sqrt{-\mathop{\rm Im} b_{-n} \mathop{\rm Im} a_m }}{\pi}\int\limits_{-\infty}^\infty
    \frac{1}{x-\overline b_{-n}}\prod_{j=1}^{-n-1}\chi_j^-\frac{x-b_j}{x-\overline b_j}\cdot
    \frac{1}{x- a_{m}}\prod_{l=0}^{m-1}\chi_l^+ \frac{x-\overline a_l}{x- a_l}~dx.
\end{equation}

The integrand in the right side of  equality (\ref{0}) is the rational fraction with poles at the points
$\overline {b}_j,~j=1,2,\ldots,-n,$ and $a_l,~l=0,1,\ldots,m,$ from the upper half plane $\mathbb{C}_+$ of complex plane. It is clear that for $|z|\to \infty,$ $ z \in \mathbb{C},$ this fraction have the order ${\cal O}(|z|^{-2}),$ that sure the convergence of the integral. Therefore, using Cauchy integral theorem for the lower half plane $\mathbb{C}_-$, we get
$$
    \frac{1}{\pi} \int\limits_{-\infty}^\infty
    \Phi_{-n}^-(x) \overline{\Phi_m^+(x)} dx=0, \quad n=-1,-2,\ldots,~m=0,1,2,\ldots.
$$

Using similar reasoning in the case $m=-1,-2,\ldots,~n= 0,1,2,\ldots,$ we completes the proof of equality (\ref{ort_Omega}).  Lemma is proved.

From formulas (\ref{Dzhrbashyan identity+}) and (\ref{Dzhrbashyan identity-}), observe that
$$
    \sum_{k=-m+1}^{n-1} \overline{\Phi_k (\zeta)} \Phi_k(z)=
    \frac{1}{2i(\overline\zeta-z)}\bigg[\prod_{k=1}^{m-1} \frac{z-b_k}{z-\overline b_k}
    \cdot\prod_{k=1}^{m-1} \frac{\overline \zeta-\overline b_k}{\overline \zeta-b_k} -
$$
$$
    -\prod_{k=0}^{n-1} \frac{z-a_k}{z-\overline a_k} \cdot
    \prod_{k=0}^{n-1} \frac{\overline \zeta-\overline a_k}{\overline \zeta-a_k}  \bigg]=
    \frac{1}{2i(\overline\zeta-z)}\Bigg[B_m^-(z) \overline{B_m^-(\zeta)}- B_n^+(z) \overline{B_n^+(\zeta)}\Bigg].
$$

\bigskip

{\bf 2. Representation of Dirichlet kernel of the system $\{\Phi_k(z)\}$ on the real axis.}
Denote the kernel of the system $\{\Phi_k(z)\}$ on the real axis $\mathbb{R}$ by
$$
    D_{n,m}({\bf a};{\bf b};x;t):=\sum_{k=-m+1}^{n-1} \overline{\Phi_k (x)} \Phi_k(t).
    \quad x,t \in \mathbb{R}.
$$

The following statement contains the representation of the quantity $D_{n,m}({\bf a};{\bf b};x;t)$ in the form, which is convenient for further research.

\begin{lem} \label{L.4}
The formula is true:
$$
    D_{n,m}({\bf a};{\bf b};x;t)=\frac{1}{2i(t-x)}
    \bigg[\mathop {\rm exp} \bigg(-2i \int\limits_x^t \sum\limits_{k=1}^{m-1}
    \frac{\mathop {\rm Im } b_k ~du}{(u- \mathop {\rm Re } b_k)^2+(\mathop {\rm Im } b_k)^2 } \bigg)-
$$
\begin{equation}\label{1}
    - \mathop {\rm exp} \bigg(-2i \int\limits_x^t \sum\limits_{k=0}^{n-1}
    \frac{\mathop {\rm Im } a_k ~du}{(u- \mathop {\rm Re } a_k)^2+(\mathop {\rm Im } a_k)^2 } \bigg)\bigg].
\end{equation}
\end{lem}

{\bf Proof.} Let $b_k=\beta_k+i\gamma_k,$ $\beta_k \in \mathbb{R}, \gamma_k <0.$
Take in account  that $x,t \in \mathbb{R},$ we find
$$
    B_m^-(x) \overline{B_m^-(t)}=\prod_{k=1}^{m-1} \frac{x-b_k}{x-\overline b_k}
    \cdot\prod_{k=1}^{m-1} \frac{ t-\overline b_k}{t-b_k}=
$$
$$
    =\prod_{k=1}^{m-1} \frac{(x-\beta_k) - i \gamma_k}{(x-\beta_k) + i \gamma_k}
    \prod_{k=1}^{m-1} \frac{(t-\beta_k) + i \gamma_k}{(t-\beta_k) - i \gamma_k}=
$$
$$
    =\prod_{k=1}^{m-1} \frac{\sqrt{(x-\beta_k)^2 + \gamma_k^2}
    (\cos \mathop{\rm arctg} \frac{\gamma_k}{x-\beta_k}-
    i\sin \mathop{\rm arctg} \frac{\gamma_k}{x-\beta_k})}{\sqrt{(x-\beta_k)^2 + \gamma_k^2}
    (\cos \mathop{\rm arctg} \frac{\gamma_k}{x-\beta_k}+
    i\sin \mathop{\rm arctg} \frac{\gamma_k}{x-\beta_k})} \times
$$
$$
    \times \prod_{k=1}^{m-1} \frac{\sqrt{(t-\beta_k)^2 + \gamma_k^2}
    (\cos \mathop{\rm arctg} \frac{\gamma_k}{t-\beta_k}+
    i\sin \mathop{\rm arctg} \frac{\gamma_k}{t-\beta_k})}{\sqrt{(t-\beta_k)^2 + \gamma_k^2}
    (\cos \mathop{\rm arctg} \frac{\gamma_k}{t-\beta_k}-
    i\sin \mathop{\rm arctg} \frac{\gamma_k}{t-\beta_k})}=
$$
$$
    =\prod_{k=1}^{m-1} \frac{(\cos \mathop{\rm arctg} \frac{\gamma_k}{x-\beta_k}-
    i\sin \mathop{\rm arctg} \frac{\gamma_k}{x-\beta_k}) (\cos \mathop{\rm arctg} \frac{\gamma_k}{t-\beta_k}+i\sin \mathop{\rm arctg} \frac{\gamma_k}{t-\beta_k})}
    {(\cos \mathop{\rm arctg} \frac{\gamma_k}{x-\beta_k}+
    i\sin \mathop{\rm arctg} \frac{\gamma_k}{x-\beta_k}) (\cos \mathop{\rm arctg} \frac{\gamma_k}{t-\beta_k}-i\sin \mathop{\rm arctg} \frac{\gamma_k}{t-\beta_k})}=
$$
$$
    =\prod_{k=1}^{m-1} \mathop{\rm exp} \bigg[ 2i\bigg(\mathop{\rm arctg} \frac{\gamma_k}{t-\beta_k}-
    \mathop{\rm arctg} \frac{\gamma_k}{x-\beta_k}   \bigg) \bigg].
$$

Note that
$$
    \bigg(\mathop{\rm arctg} \frac{\gamma_k}{\cdot-\beta_k}\bigg)'=
    -\frac{\gamma_k}{(\cdot-\beta_k)^2+\gamma_k^2},
$$
we obtain
$$
    \prod_{k=1}^{m-1} \frac{x-b_k}{x-\overline b_k}
    \cdot\prod_{k=1}^{m-1} \frac{ t-\overline b_k}{t-b_k}=
    \prod_{k=1}^{m-1} \mathop{\rm exp} \bigg[ 2i\bigg(\mathop{\rm arctg} \frac{\gamma_k}{t-\beta_k}-
    \mathop{\rm arctg} \frac{\gamma_k}{x-\beta_k}   \bigg) \bigg]=
$$
$$
    =\mathop{\rm exp} \bigg[ 2i \sum\limits_{k=1}^{m-1} \bigg(\mathop{\rm arctg} \frac{\gamma_k}{t-\beta_k}-
    \mathop{\rm arctg} \frac{\gamma_k}{x-\beta_k}   \bigg) \bigg]=
$$
$$
    =\mathop{\rm exp} \bigg[ -2i \sum\limits_{k=1}^{m-1} \int\limits_x^t
    \frac{\gamma_k du}{(u-\beta_k)^2+\gamma_k^2} \bigg].
$$

Using similar reasoning for the value
$
    B_n^+(z) \overline{B_n^+(\zeta)},
$
we get (\ref{1}). Lemma is proved.

\bigskip

{\bf 3. Statement of the problem and historical review.}
Let, as above,
$$
    {\bf a}:=\{a_k\}_{k=0}^{\infty}, ~ (\mathop{\rm Im} a_k>0), \quad
    {\bf b}:=\{b_k\}_{k=1}^{\infty}, ~ (\mathop{\rm Im} b_k<0)
$$
be arbitrary sequences of complex numbers from upper  $\mathbb C_+$ and lower
$\mathbb C_-$ half plane of the complex plane $\mathbb C$ respectively.

In the paper \cite{Kober_1944} it was shown that the system of the rational functions
\begin{equation}\label{element-racfunct}
    \frac{1}{z-\bar{a}_0}, \quad \frac{1}{z-\bar{b}_1}, \quad \frac{1}{z-\bar{a}_1}, \quad \frac{1}{z-\bar{b}_1}, \ldots,
\end{equation}
is closed with respect to $L_p(\mathbb{R})~ (1<p<\infty)$ if, and only if, the series
$$
    \sigma({\bf a}):=\sum\limits_{k=0}^\infty \frac{|\mathop {\rm Im } a_k|}{1+|a_k|^2}, \quad
    \sigma({\bf b}):=\sum\limits_{k=0}^\infty \frac{|\mathop {\rm Im } b_k|}{1+|b_k|^2},
$$
diverge.  The orthogonalization of the system (\ref{element-racfunct}) on the real axis $\mathbb{R}$ leads to the system $\{\Phi_n(x)\},~x \in \mathbb{R}, n \in \mathbb{Z},$ therefore to arbitrary function from $f \in L_2(\mathbb{R})$ can be put  to conformity its Fourier series on the system $\{\Phi_n(x)\}$:
\begin{equation}\label{Furier-series}
    f(x) \sim \sum\limits_{k=-\infty}^\infty c_k \Phi_k (x), \quad x \in \mathbb{R},
\end{equation}
where
$$
    c_k=\frac{1}{\pi} \int\limits_{-\infty}^\infty f(x) \overline{\Phi_k(x)}dx,
    \quad k=0,\pm 1, \pm 2, \ldots,
$$
whose partial sums
$$
    S_{n,m}(f;{\bf a};{\bf b};x)=\sum\limits_{k=-m+1}^{n-1} c_k \Phi_k(x),
$$
converge in mean square to the functions $f(x),$   i.e.
$$
    \lim\limits_{n,m\to\infty} \frac{1}{\pi} \int\limits_{-\infty}^\infty |f(x)-S_{n,m}(f;x)|^2dx=0,
$$
on condition that  series $\sigma({\bf a})$ and $\sigma({\bf b})$ diverge.

However, formal Fourier series of form (\ref{Furier-series}) can be written for any function
$f \in L_1 (\mathbb{R})$ and for partial sums of this series will have
$$
    S_{n,m} (f;{\bf a};{\bf b};x)=\sum\limits_{k=-m}^n c_k \Phi_k(x)=
    \frac{1}{\pi} \int\limits_{-\infty}^\infty f(t) \bigg[
    \sum\limits_{k=-m}^n \overline{\Phi_k(t)} \Phi_k(x) \bigg] dt=
$$
\begin{equation}\label{represent}
    =\frac{1}{\pi} \int\limits_{-\infty}^\infty f(t) D_{n,m}({\bf a};{\bf b};x;t)~dt.
\end{equation}

Our aim is to investigate the problems of convergence in the metrics of the spaces $L_p(\mathbb{R})$ $(1<p<\infty)$ and pointwise convergence of the partial sums $S_{n,m} (f;{\bf a};{\bf b};x)$ of  Fourier series (\ref{Furier-series}) to the corresponding function $f$ as $m,n \to \infty.$

Note that the extremal problems of the approximation on the real axis by rational functions with fixed poles was originated from the works of S.N.~Bernshtain. So in the monograph \cite{Bernshtain_Extrem_properties} was built the rational fraction which deviates least from zero on the real axis, and also in the first time was obtained the solution of the problem of approximation of functions by rational fractions in the uniform metric.  Researches of  conditions of completeness of the system of  the functions (\ref{element-racfunct}), conditions of convergence of the approximating aggregate by these systems in the metrics of the spaces $L_p(\mathbb{R}),~p\in[1;\infty],$   were considered in \cite{Kober_1944}.  N.I. Achieser \cite{Ahiezer.lect.aprox},
exploring weighted polynomial approximation problem, first establish the exact value of the best weighted approximation on the real axis the kernels of form
\begin{equation}\label{kernel-2}
    \frac{Ax+D}{x^2+\lambda}, \quad (\mathop {\rm Im } A=\mathop {\rm Im } D=0, ~ \lambda>0)
\end{equation}
for the weight of a certain kind.

For the first time, the Fourier series expansion on the system (\ref{Phi-system}) was studied by
M.M.~Dzhrbashyan \cite{Djrbashan_1974}. He developed the method that allowed to obtain the solving of extreme problems of best rational approximation of the Cauchy kernel
$$
\frac{1}{\zeta-z}, \quad \mathop {\rm Im } \zeta\not=0, \quad z\in \mathbb{R},
$$
both in uniform metric and in the mean square metric. The proposed method is based on the use of orthogonal system (\ref{Phi-system}) and certain biorthogonal system of rational functions with fixed poles on the real axis $\mathbb{R}$. Using the method of M.M.~Dzhrbashyan in the paper \cite{Voskanian}  was solved similar problems for the kernels of kind (\ref{kernel-2}). In the works of  V.M.~Rusak (see monograph \cite{Rusak_Racion_funct}) were built the rational operators of type Fejer, Valle Poussin and Jackson and investigated approximation properties of these operators. In the article \cite{Savchuk_CH_umg_2014} computed the value of the best approximation of the Cauchy kernel on the real axis $\mathbb R$ by some subspaces from $L_q(\mathbb R)$. This result is applied to the evaluation of the exact upper bounds for pointwise deviation of certain interpolation operators with interpolation nodes in the upper half plane and certain linear means of Fourier series on the Takenaka-Malmquist system from the functions in a unit ball of the Hardy space $H_p$, $2\le p<\infty$.

For simplicity reasons, in this paper we assume that the relevant terms of the sequences ${\bf a}=\{a_k\}_{k=0}^{\infty}$ and ${\bf b}=\{b_k\}_{k=1}^{\infty}$ are  pairwise  conjugated, i.e.
$b_k=\bar{a}_{k-1},$ $(k=1,2,\ldots)$. Then, the representation of the spectral function (\ref{1}) and the partial sums (\ref{represent}) in case when $m=n+1$ take the form
\begin{equation}\label{2}
    D_{n,n+1}({\bf a};{\bf \bar a};x;t)=\frac{1}{(t-x)} \sin \bigg(\int\limits_x^t \bigg[ \sum\limits_{k=0}^{n-1}
    \frac{2\mathop {\rm Im } a_k}{(u- \mathop {\rm Re } a_k)^2+(\mathop {\rm Im } a_k)^2 } \bigg]du \bigg),
\end{equation}
and
$$
    S_{n,n+1} (f;{\bf a};{\bf \bar a};x):=S_{n} (f;{\bf a};x)=\frac{1}{\pi} \int\limits_{-\infty}^\infty f(t) D_{n,n+1} ({\bf a}; {\bf \bar a}; x;t)~dt,
$$
where ${\bf a}=\{a_k\}_{k=0}^{\infty}$ and ${\bf \bar a}=\{\bar{a}_{k-1}\}_{k=1}^{\infty}.$

\bigskip

{\bf 4.  Main results.} Now we can formulate the main results of this paper.

\begin{teo} \label{T.1}
Suppose that sequence ${\bf a}:=\{a_k\}_{k=0}^{\infty}$ satisfies the condition
\begin{equation}\label{condition-Blyashke}
    \lim_{n\to\infty} \sum_{k=0}^{n-1} \frac{|\mathop {\rm Im } a_k|}{1+|a_k|^2} =+\infty.
\end{equation}

Then for any function $f \in L_p(\mathbb{R}),~1<p<\infty,$ its Fourier series on the system (\ref{system}) converges to this function in the metric of the spaces $L_p(\mathbb{R}),$ i.e.:
$$
    \lim_{n\to \infty} \int\limits_{-\infty}^\infty |f(x)-S_{n} (f;{\bf a};x)|^p dx=0,
    \quad 1<p<\infty.
$$
\end{teo}

Let
\begin{equation}\label{sigma}
    \sigma_n:=\sum_{k=0}^{n-1} \frac{|\mathop {\rm Im } a_k|}{1+|a_k|^2}, \quad
    \varsigma_n:=\sum_{k=0}^{n-1} \frac{1}{(\mathop {\rm Im } a_k)^2}.
\end{equation}

\begin{teo} \label{T.2}
 Assume that function $f \in L_1(\mathbb{R})$ and has bounded variation on $\mathbb{R}$. If the sequence ${\bf a}:=\{a_k\}_{k=0}^{\infty}$ has no limit points on the real axis $\mathbb{R},$ satisfies the condition (\ref{condition-Blyashke}) and $\varsigma_n/\sigma_n\le {\rm const}$, then at each point $x_0 \in \mathbb{R}$ the equality is true:
$$
    \lim\limits_{n\to \infty} S_{n}(f;{\bf a};x_0)= \frac{f(x_0-0)+f(x_0+0)}{2}.
$$
\end{teo}

\begin{cor}\label{N.2}
If all conditions of the theorem \ref{T.2} are satisfied and $x_0$ is the point of continuity of function $f,$ then
$$
    \lim\limits_{n\to \infty} S_{n}(f;{\bf a};x_0)=f(x_0).
$$
\end{cor}

\begin{teo} \label{T.3}
Let function $f \in L_1(\mathbb{R})$ and the limit values $f(x_0-0)$ and $f(x_0+0)$ exist at the point $x_0$. If the integrals
\begin{equation}\label{T3-cond}
    \int\limits_0^\delta \frac{f(x_0-y)-f(x_0-0)}{y}dy, \quad
    \int\limits_0^\delta \frac{f(x_0+y)-f(x_0+0)}{y}dy,
\end{equation}
exist and the sequence ${\bf a}:=\{a_k\}_{k=0}^{\infty}$ has no limit points on the real axis $\mathbb{R}$, satisfies the condition (\ref{condition-Blyashke}) and $\varsigma_n/\sigma_n\le {\rm const}$, then
$$
    \lim\limits_{n\to \infty} S_{n}(f;{\bf a};x_0)= \frac{f(x_0-0)+f(x_0+0)}{2}.
$$
\end{teo}

It is obvious that sequence $ {\bf a}: = \{a_k \}_{k = 0}^{\infty} $  satisfies the conditions of Theorems \ref{T.2} and \ref{T.3}, if we can find positive constants $ C_1 $ and $ C_2$, which do not depend on $k$ and such that $ 0 <C_1 \le | a_k | \le C_2, ~ k \in \mathbb {N}$.
Note that the conditions of these Theorems can be performed also in the case, when  $ | a_k | \to \infty, ~ n \to \infty $. For example, if the real and imaginary parts of $ {\bf a}: = \{a_k \}_{k = 0}^{\infty} $  satisfy the relations
$$
    |\mathop {\rm Re } a_k|={\cal O}(k^\alpha), ~ 0\le\alpha \le {3\over4}, \quad
    |\mathop {\rm Im } a_k|={\cal O}(k^\beta),~ {1\over2} < \beta\le 1, \quad k \to \infty,
$$
it is easy to verify, that this sequence  has no limit points on the real axis $ \mathbb {R} $, satisfies the condition (\ref{condition-Blyashke}) and $ \varsigma_n / \sigma_n \le {\rm const}.$

{\bf 5. Auxiliary results.} We first obtain following statement, where  the integral representation for the partial sums $S_{n}(f;{\bf a};x)$ of Fourier series on the system (\ref{system}) has convenient form.

\begin{lem}\label{L.5}
If $f \in L_1(\mathbb{R})$, then for any $x \in \mathbb{R}$ the equality holds
\begin{equation}\label{Sn}
    S_{n}(f;{\bf a};x)=\frac{1}{\pi} \int\limits_0^\infty f(x-y)
    \frac{\sin y \mu_n(-y;x)}{y}~dy+
$$
$$
    +\frac{1}{\pi} \int\limits_0^\infty f(x+y)
    \frac{\sin y \mu_n(y;x)}{y}~dy,
\end{equation}
where
\begin{equation}\label{def-mu}
    \mu_n(y;x):=\frac{1}{y} \int\limits_x^{x+y} \bigg[\sum\limits_{k=0}^{n-1}
    \frac{2\mathop {\rm Im } a_k}{(u- \mathop {\rm Re } a_k)^2+(\mathop {\rm Im } a_k)^2 } \bigg]du.
\end{equation}
\end{lem}

{\bf Proof.} We split the integral on the whole axis to two parts over intervals $(-\infty;x)$ and $(x;\infty),$ and  make the change of variables. Given the relations (\ref{2}),  we find
$$
    S_{n}(f;{\bf a};x)=\frac{1}{\pi} \int\limits_{-\infty}^\infty f(t) D_{n,n+1}({\bf a};{\bf \bar a};t;x)~dt=
$$
$$
    =\frac{1}{\pi} (\int\limits_{-\infty}^x+\int\limits_{x}^\infty) f(t) D_{n,n+1}({\bf a};{\bf \bar a};t;x)~dt=
$$
$$
    =\frac{1}{\pi} \int\limits_{0}^\infty f(x-y) D_{n,n+1}({\bf a};{\bf \bar a};x-y;x)~dy+
$$
$$
    +\frac{1}{\pi} \int\limits_{0}^\infty f(x+y) D_{n,n+1}({\bf a};{\bf \bar a};x+y;x)~dy=
$$
$$
    =\frac{1}{\pi} \int\limits_{0}^\infty \frac{f(x-y)}{y}
    \sin \bigg( \int\limits_{x-y}^{x}\bigg[\sum\limits_{k=0}^{n-1}
    \frac{2\mathop {\rm Im } a_k}
    {(u- \mathop {\rm Re } a_k)^2+(\mathop {\rm Im } a_k)^2 } \bigg]du \bigg)dy+
$$
$$
    +\frac{1}{\pi} \int\limits_{0}^\infty \frac{f(x+y)}{y}
    \sin \bigg(\int\limits_{x}^{x+y}\bigg[\sum\limits_{k=0}^{n-1}
    \frac{2\mathop {\rm Im } a_k}
    {(u- \mathop {\rm Re } a_k)^2+(\mathop {\rm Im } a_k)^2 } \bigg]du \bigg)dy.
$$

Hence, considering the designation (\ref{def-mu}), we obtain the equality (\ref{Sn}). Lemma is proved.

\begin{lem}\label{L.6}
For each fixed $x\in \mathbb{R}$ and $y>0$ the inequalities hold:
\begin{equation}\label{estim-ymu'}
    \bigg|[y\mu_n(\pm y; x)]'_y \bigg|\ge \frac{1}{1+(|x|+y)^2} \sum_{k=0}^{n-1}
     \frac{|\mathop {\rm Im } a_k|}{1+|a_k|^2},
\end{equation}
\begin{equation}\label{estim-mu}
     |\mu_n(\pm y; x)| \ge \frac{1}{1+(|x|+y)^2} \sum_{k=0}^{n-1}
     \frac{|\mathop {\rm Im } a_k|}{1+|a_k|^2},
\end{equation}

\end{lem}

{\bf Proof.}  On the basis of the relation (\ref{def-mu}), we find
$$
    [y\mu_n(y; x)]'_y= \frac{d}{dy}\int\limits_x^{x+y} \bigg[\sum\limits_{k=0}^{n-1}
    \frac{2\mathop {\rm Im } a_k}{(u -\mathop {\rm Re } a_k)^2+(\mathop {\rm Im } a_k)^2 } \bigg]du =
$$
\begin{equation}\label{3}
    = \sum\limits_{k=0}^{n-1}
    \frac{2\mathop {\rm Im } a_k}{([y+x]- \mathop {\rm Re } a_k)^2+(\mathop {\rm Im } a_k)^2 }.
\end{equation}
Hence, given the inequality
\begin{equation}\label{estim-low}
    (1+(|x|+y)^2)(1+(\mathop {\rm Re } a_k)^2+(\mathop {\rm Im } a_k)^2)\ge
    ([t+x]-\mathop {\rm Re } a_k)^2+(\mathop {\rm Im } a_k)^2,
\end{equation}
which holds for any $x\in \mathbb{R}$, $y>0$  and $0 \le t \le y,$ we get the estimation  (\ref{estim-ymu'}) for the value $|[y\mu_n(y; x)]'_y|.$

Since
$$
    \mu_n(y; x)=\frac{1}{y} \int\limits_x^{x+y} \bigg[\sum\limits_{k=0}^{n-1}
    \frac{2\mathop {\rm Im } a_k}{(u- \mathop {\rm Re } a_k)^2+(\mathop {\rm Im } a_k)^2 } \bigg]du=
$$
$$
    =\frac{1}{y} \int\limits_0^{y} \bigg[\sum\limits_{k=0}^{n-1}
    \frac{2\mathop {\rm Im } a_k}{([t+x]- \mathop {\rm Re } a_k)^2+(\mathop {\rm Im } a_k)^2 } \bigg]dt,
$$
then  applying the inequality (\ref{estim-low}) again, we
obtain the inequality  (\ref{estim-mu}) for the value $|\mu_n(y; x)|.$

The values $|[y\mu_n(-y; x)]'_y|$ and $\mu_n(-y; x)$ are estimated similarly. Lemma is proved.

\begin{lem}\label{L.6'} Uniformly for $x \in \mathbb{R}$ and $y>0$ the estimations are true:
\begin{equation}\label{estim-ymu''}
    \bigg|[y\mu_n(\pm y; x)]''_{y^2} \bigg|\le \sum_{k=0}^{n-1} \frac{1}{(\mathop {\rm Im } a_k)^2},
\end{equation}
\begin{equation}\label{estim-mu'}
    \bigg|[\mu_n(\pm y; x)]'_y \bigg|\le \sum_{k=0}^{n-1} \frac{1}{(\mathop {\rm Im } a_k)^2},
\end{equation}
\begin{equation}\label{estim-mu''}
    \bigg|[\mu_n(\pm y; x)]''_{y^2} \bigg|\le \frac{8}{3} \sum_{k=0}^{n-1}
    \frac{1}{(\mathop {\rm Im } a_k)^3}.
\end{equation}
\end{lem}

{\bf Proof.}  On the basis of the relation (\ref{3}), we find
$$
    [y\mu_n(y; x)]''_{y^2}= \frac{d}{dy} \Bigg[\sum\limits_{k=0}^{n-1}
    \frac{2\mathop {\rm Im } a_k}{([x+y]- \mathop {\rm Re } a_k)^2+(\mathop {\rm Im } a_k)^2 }\Bigg] =
$$
$$
    = -\sum\limits_{k=0}^{n-1}
    \frac{4\mathop {\rm Im } a_k ([y+x]- \mathop {\rm Re } a_k)}
    {\bigg(([y+x]- \mathop {\rm Re } a_k)^2+(\mathop {\rm Im } a_k)^2\bigg)^2 }.
$$
Hence, taking into account the obvious inequalities
\begin{equation}\label{estim-up1}
    \frac{2\mathop {\rm Im } a_k(v- \mathop {\rm Re } a_k)}
    {(v- \mathop {\rm Re } a_k)^2+(\mathop {\rm Im } a_k)^2} \le 1,
\end{equation}
and
\begin{equation}\label{estim-up}
    \frac{1} {(v- \mathop {\rm Re } a_k)^2+(\mathop {\rm Im } a_k)^2 } \le
    \frac{1}{(\mathop {\rm Im } a_k)^2} ,
\end{equation}
at $v=t+x$, we obtain the estimation (\ref{estim-ymu''}) for the value $|[y\mu_n(y; x)]''_{y^2}|$.

To prove the estimation (\ref{estim-mu'}) it should be noted that
$$
    [\mu_n(y; x)]'_y=\frac{1}{y^2} \int\limits_0^y \bigg[\sum\limits_{k=0}^{n-1}
    \frac{2\mathop {\rm Im } a_k}{([y+x]- \mathop {\rm Re } a_k)^2+(\mathop {\rm Im } a_k)^2 }-
$$
$$
    -\sum\limits_{k=0}^{n-1} \frac{2\mathop {\rm Im } a_k}{([t+x]- \mathop {\rm Re } a_k)^2+
    (\mathop {\rm Im } a_k)^2 } \bigg]dt=
$$
$$
    =\frac{1}{y^2} \int\limits_0^y \sum\limits_{k=0}^{n-1}
    \frac{2\mathop {\rm Im } a_k(t-y) [(y+x- \mathop {\rm Re } a_k)+(t+x- \mathop {\rm Re } a_k)]}
    {[(y+x- \mathop {\rm Re } a_k)^2+(\mathop {\rm Im } a_k)^2]
    [(t+x- \mathop {\rm Re } a_k)^2+(\mathop {\rm Im } a_k)^2]}=
$$
$$
    =\frac{1}{y^2} \int\limits_0^y (t-y) \sum\limits_{k=0}^{n-1} \bigg[
    \frac{2\mathop {\rm Im } a_k(y+x- \mathop {\rm Re } a_k)}
    {(y+x- \mathop {\rm Re } a_k)^2+(\mathop {\rm Im } a_k)^2} \cdot
     \frac{1}{(t+x- \mathop {\rm Re } a_k)^2+(\mathop {\rm Im } a_k)^2}+
$$
$$
    +\frac{2\mathop {\rm Im } a_k(t+x- \mathop {\rm Re } a_k)}
    {(t+x- \mathop {\rm Re } a_k)^2+(\mathop {\rm Im } a_k)^2} \cdot
     \frac{1}{(y+x- \mathop {\rm Re } a_k)^2+(\mathop {\rm Im } a_k)^2} \bigg]dt.
$$
Hence, again using the inequalities (\ref{estim-up1}) -- (\ref{estim-up})
at $v=t+x$ and $v=y+x$, we get
$$
    \bigg|[\mu_n(y; x)]'_y\bigg| \le 2 \sum_{k=0}^{n-1} \frac{1}{(\mathop {\rm Im } a_k)^2}\cdot
     \frac{ 1}{y^2} \int\limits_0^y (y-t)dt=\sum_{k=0}^{n-1} \frac{1}{(\mathop {\rm Im } a_k)^2},
$$
which proves the estimation (\ref{estim-mu'}) for the value $|[\mu_n(y; x)]'_y|$.

Then we find
$$
    [\mu_n(y;x)]''_{y^2}=\frac{d}{dy}\Bigg[\frac{1}{y} \sum\limits_{k=0}^{n-1}
    \frac{2\mathop {\rm Im } a_k}{([y+x]- \mathop {\rm Re } a_k)^2+(\mathop {\rm Im } a_k)^2 }-
$$
$$
    -\frac{1}{y^2}\int\limits_0^y \sum\limits_{k=0}^{n-1} \frac{2\mathop {\rm Im } a_k ~dt}{([t+x]- \mathop {\rm Re } a_k)^2+
    (\mathop {\rm Im } a_k)^2 } \Bigg]=
$$
$$
    =\frac{2}{y^3}\int\limits_0^y \sum\limits_{k=0}^{n-1} \frac{2\mathop {\rm Im } a_k ~dt}{([t+x]- \mathop {\rm Re } a_k)^2+
    (\mathop {\rm Im } a_k)^2 } -
$$
$$
    -\frac{2}{y^2} \sum\limits_{k=0}^{n-1}
    \frac{2\mathop {\rm Im } a_k}{([y+x]- \mathop {\rm Re } a_k)^2+
    (\mathop {\rm Im } a_k)^2 }
$$
$$
    -\frac{1}{y} \sum\limits_{k=0}^{n-1} \frac{4\mathop {\rm Im } a_k (x+y-\mathop {\rm Re } a_k)}
    {([t+x]- \mathop {\rm Re } a_k)^2+(\mathop {\rm Im } a_k)^2 },
$$
where after elementary transformations, we obtain
$$
    [\mu_n(y;x)]''_{y^2}=\frac{2}{y^3} \int\limits_0^y (y-t) \sum\limits_{k=0}^{n-1}
    \frac{2\mathop {\rm Im } a_k}
    {(y+x- \mathop {\rm Re } a_k)^2+(\mathop {\rm Im } a_k)^2 } \times
$$
$$
    \times\Bigg[ \frac{y+t+2x-2 \mathop {\rm Re } a_k}{(t+x- \mathop {\rm Re } a_k)^2+
    (\mathop {\rm Im } a_k)^2 } -\frac{2(y+x- \mathop {\rm Re } a_k)}{([y+x]- \mathop {\rm Re } a_k)^2+
    (\mathop {\rm Im } a_k)^2 }\Bigg]dt=
$$
$$
    =\frac{2}{y^3} \int\limits_0^y \sum\limits_{k=0}^{n-1}
    \Bigg[ \frac{4\mathop {\rm Im } a_k(y+x- \mathop {\rm Re } a_k)^2}
    {[(t+x- \mathop {\rm Re } a_k)^2+
    (\mathop {\rm Im } a_k)^2 ] [(y+x- \mathop {\rm Re } a_k)^2 +(\mathop {\rm Im } a_k)^2]^2 }+
$$
$$
    +\frac{4\mathop {\rm Im } a_k(y+x- \mathop {\rm Re } a_k)(t+x- \mathop {\rm Re } a_k)}
    { [(t+x- \mathop {\rm Re } a_k)^2+
    (\mathop {\rm Im } a_k)^2 ][(y+x- \mathop {\rm Re } a_k)^2 +(\mathop {\rm Im } a_k)^2]^2}-
$$
$$
    -\frac{2\mathop {\rm Im } a_k}{[(y+x- \mathop {\rm Re } a_k)^2 +(\mathop {\rm Im } a_k)^2]
    [(t+x- \mathop {\rm Re } a_k)^2+(\mathop {\rm Im } a_k)^2]}\Bigg](y-t)^2dt.
$$
Finally, applying the inequality (\ref{estim-up1}) -- (\ref{estim-up})
at $v=t+x$ and $v=y+x$, we get
$$
    |[\mu_n(y;x)]''_{y^2}|\le \frac{2}{y^3} \int\limits_0^y  \sum\limits_{k=0}^{n-1}
    \frac{4(y-t)^2}{(\mathop {\rm Im } a_k)^3}~dt = \frac{8}{3} \sum\limits_{k=0}^{n-1}
    \frac{1}{(\mathop {\rm Im } a_k)^3},
$$
which proves the inequality (\ref{estim-mu''}) for the value $|[\mu_n(y; x)]''_{y^2}|$. The values $|[y\mu_n(-y; x)]''_{y^2}|$, $|[\mu_n(-y; x)]'_y|$ and $|[\mu_n(-y; x)]''_{y^2}|$ are estimated similarly. Lemma is proved.

\begin{lem}\label{L.7'}
Suppose that sequence ${\bf a}:=\{a_k\}_{k=0}^{\infty}$ has no limit points on the real axis $\mathbb{R}$, satisfies the condition (\ref{condition-Blyashke}) and $\varsigma_n/\sigma_n\le {\rm const}$. Then, for arbitrary  function $\varphi \in L(\mathbb{R}_+),$ the equation holds
\begin{equation}\label{rest}
    \lim_{n \to \infty} \int\limits_0^{\infty} \varphi(y) \sin [y \mu_n(\pm y; x)]dy=0.
\end{equation}
\end{lem}

{\bf Proof.} We fix number $\varepsilon>0.$ Obviously, for arbitrary function $\varphi \in L(\mathbb{R}_+),$  the number $\theta=\theta(\varepsilon)>0$ such that
\begin{equation}\label{rest-1}
    \int\limits_\theta^{\infty} |\varphi(y)| dy \le \frac{\varepsilon}{3}
\end{equation}
exists always.

By K.~Weierstrass theorem on the interval $[0; \theta]$ can be found polynomial $P(y)=P(\varphi;y),$ which provides the estimate
\begin{equation}\label{rest-2}
    \int\limits_0^\theta |\varphi(y)- P(y)| dy \le \frac{\varepsilon}{3}.
\end{equation}

From the relations (\ref{rest-1}) and (\ref{rest-2}) it follows
$$
    \Bigg|\int\limits_0^{\infty} \varphi(y) \sin [y \mu_n(\pm y; x)]dy \Bigg| = \Bigg|
    \int\limits_0^\theta  P(y) \sin [y \mu_n(\pm y; x)]dy+
$$
$$
    +\int\limits_\theta^{\infty} \varphi(y) \sin [y \mu_n(\pm y; x)]dy+
    \int\limits_0^\theta  [\varphi(y)-P(y)] \sin [y \mu_n(\pm y; x)]dy\Bigg|\le
$$
$$
    \le \Bigg| \int\limits_0^\theta  P(y) \sin [y \mu_n(\pm y; x)]dy \Bigg|+\frac{2 \varepsilon}{3},
$$
herewith uniformly over $n=0,1,2,\ldots$.

So to complete the proof of Lemma left to show that we can find a positive integer $n_0=n_0(\varepsilon)>0$, such that for all natural $n > n_0$ the inequality true
\begin{equation}\label{rest-3}
    \Bigg| \int\limits_0^\theta  P(y) \sin [y \mu_n(\pm y; x)] dy \Bigg| \le \frac{\varepsilon}{3}.
\end{equation}

Let
$$
    M_1= \max_{0\le y \le \theta} |P(y)|, \quad
    M_2= \max_{0\le y \le \theta} |P'(y)|.
$$

Integrating by parts, we find
$$
    \int\limits_0^\theta  P(y) \sin [y \mu_n(\pm y; x)] dy= \left.
    -\frac{P(y) }{[y \mu_n(\pm y; x)]'_y}\cos [y \mu_n(\pm y; x)] \right|_0^\theta+
$$
$$
    +\int\limits_0^\theta \frac{P'(y) }{[y \mu_n(\pm y; x)]'_y}\cos [y \mu_n(\pm y; x)]dy-
$$
$$
    -\int\limits_0^\theta \frac{P(y) [y \mu_n(\pm y; x)]''_{y^2}}{([y \mu_n(\pm y; x)]'_y)^2}
    \cos [y \mu_n(\pm y; x)]dy.
$$

We apply the estimations (\ref{estim-ymu'}) -- (\ref{estim-ymu''}) for integrands in the right side of  last equality. Considering the designation (\ref{sigma}), we get
$$
    \Bigg|\int\limits_0^\theta  P(y) \sin [y \mu_n(\pm y; x)] dy \Bigg| \le
    \frac{2M_1(1+(|x|+\theta)^2)}{\sigma_n}+
$$
$$
    +\frac{M_2\theta (1+(|x|+\theta)^2)}{\sigma_n}+
    \frac{M_1 \theta (1+(|x|+\theta)^2)^2 \varsigma_n}{\sigma_n^2}.
$$

Hence, in view of (\ref{condition-Blyashke}) and the conditions $\varsigma_n/\sigma_n\le {\rm const}$ we obtain the relation (\ref{rest-3}). Lemma is proved.

\begin{lem}\label{L.7}
Assume that sequence ${\bf a}:=\{a_k\}_{k=0}^{\infty}$ has no limit points on the real axis $\mathbb{R}$, satisfies the condition (\ref{condition-Blyashke}) and $\varsigma_n/\sigma_n\le {\rm const}$. Then for arbitrary number $\delta>0$ at each fixed $x \in \mathbb{R}$
\begin{equation}\label{pi:2}
    \lim_{n \to \infty} \int\limits_0^{\delta} \frac{\sin [y \mu_n(\pm y; x)]}{y}dy=\frac{\pi}{2}.
\end{equation}
\end{lem}

{\bf Proof.} We fix $x \in \mathbb{R}$ and number $\delta>0.$  We have the equality
\begin{equation}\label{I(nx)}
    \int\limits_0^{\delta} \frac{\sin [y \mu_n(\pm y; x)]}{y}dy =
    \int\limits_0^{\delta} \frac{\sin [y \mu_n(\pm y; x)]}{y \mu_n(\pm y; x)}d[y \mu_n(\pm y; x)]-
$$
$$
    -\int\limits_0^{\delta} \sin [y \mu_n(\pm y; x)] \frac{[\mu_n(\pm y; x)]'}{ \mu_n(\pm y; x)}dy:=I_1(n;x)-I_2(n;x).
\end{equation}

Denoting  $v=y \mu_n(\pm y; x),$ we get
\begin{equation}\label{I1}
    I_1(n;x)=\int\limits_0^{\delta \mu_n(\pm\delta;x)} \frac{\sin v}{v}dv.
\end{equation}
Using the inequality (\ref{estim-low}), we find
$$
    y\mu_n(\pm y;x) =\int\limits_0^{y} \bigg[\sum\limits_{k=0}^{n-1}
    \frac{2\mathop {\rm Im } a_k}{([t\pm x]- \mathop {\rm Re } a_k)^2+(\mathop {\rm Im } a_k)^2 } \bigg]dt \ge
$$
$$
    \ge \frac{y}{1+(|x|+y)^2} \sum_{k=0}^{n-1}
     \frac{|\mathop {\rm Im } a_k|}{1+|a_k|^2}, \quad y>0,
$$
from which in view of the condition (\ref{condition-Blyashke}) implies, that  the relation is true
$$
    \delta \mu_n(\pm \delta; x)\ge \frac{\delta}{1+(|x|+\delta)^2} \sum_{k=0}^{n-1}
     \frac{|\mathop {\rm Im } a_k|}{1+|a_k|^2}\to \infty, \quad n \to \infty.
$$

Taking into account this fact, we have
\begin{equation}\label{limI1}
    \lim\limits_{n\to \infty} I_1(n;x)=\lim\limits_{n\to \infty}\int\limits_0^{\delta \mu_n(\pm\delta;x)} \frac{\sin v}{v}dv=
    \int\limits_0^{\infty} \frac{\sin v}{v}dv=\frac{\pi}{2}.
\end{equation}

Further, integrating by parts, we obtain
$$
    |I_2(n;x)|=\left.\frac{ [\mu_n(\pm y; x)]'_y\cos [y \mu_n(\pm y; x)]}
    {\mu_n(\pm y; x) [y \mu_n(\pm y; x)]'_y}\right|_0^\delta-
$$
$$
    -\int\limits_0^{\delta} \cos [y \mu_n(\pm y; x)] \frac{[\mu_n(\pm y; x)]''_{y^2}}
    { \mu_n(\pm y; x) [y \mu_n(\pm y; x)]'_y}dy+
$$
$$
    +\int\limits_0^{\delta} \cos [y \mu_n(\pm y; x)] \frac{\bigg([\mu_n(\pm y; x)]'_{y}\bigg)^2}
    { [\mu_n(\pm y; x)]^2 [y \mu_n(\pm y; x)]'_y}dy+
$$
\begin{equation}\label{4}
    +\int\limits_0^{\delta} \cos [y \mu_n(\pm y; x)] \frac{[\mu_n(\pm y; x)]'_{y}[y \mu_n(\pm y; x)]''_{y^2}}{ \mu_n(\pm y; x) ([y \mu_n(\pm y; x)]'_y)^2}dy.
\end{equation}

By the condition of Theorem the sequence ${\bf a}:=\{a_k\}_{k=0}^{\infty}$ has no limit points on the real axis $\mathbb{R}$, therefore
$$
    \sum\limits_{k=0}^{n-1}
    \frac{1}{(\mathop {\rm Im } a_k)^3}\le C \varsigma_n.
$$
Considering this fact and applying the inequalities (\ref{estim-ymu'}) -- (\ref{estim-mu}) and also (\ref{estim-ymu''}) -- (\ref{estim-mu''}) for integrands in the right side of (\ref{4}), we find
\begin{equation}\label{I2}
    |I_2(n;x)| \le C\Bigg( \frac{\varsigma_n}{\sigma_n^2} [1+(|x|+\delta)^2] +
    \frac{\varsigma_n}{\sigma_n^2} [1+(|x|+\delta)^2]^2 \delta+
    \frac{\varsigma_n^2}{\sigma_n^3} [1+(|x|+\delta)^2]^3\delta\Bigg),
\end{equation}
where $\varsigma_n$ and $\sigma_n$ are sequences, which defined by relations (\ref{sigma}).

Since $\varsigma_n/\sigma_n\le {\rm const}$, then at each fixed $x \in \mathbb{R}$ will have
\begin{equation}\label{limI2}
    \lim_{n\to \infty}|I_2(n;x)| =0.
\end{equation}

Combining the equality  (\ref{I(nx)}),   (\ref{limI1}) and (\ref{limI2}), we obtain the assertion of lemma. Lemma is proved.

Finally we prove the following lemma.

\begin{lem}\label{L.8}
Suppose that sequence ${\bf a}:=\{a_k\}_{k=0}^{\infty}$ satisfies the condition (\ref{condition-Blyashke}) and $\varsigma_n/\sigma_n\le {\rm const}$.  If the function $g(y) \in L(\mathbb{R}_+)$ and increases monotonically on the interval $[0;\infty),$ then for an arbitrary number $\delta>0$ at each fixed $x \in \mathbb{R}$
\begin{equation}\label{gpi:2}
    \lim_{n \to \infty} \int\limits_0^{\delta} g(y)\frac{\sin [y \mu_n(\pm y; x)]}{y}dy=
    \frac{\pi}{2}g(+0).
\end{equation}

\end{lem}

{\bf Proof.} We fix $x \in \mathbb{R}$ and number $\delta>0.$ Since
$$
    \int\limits_0^{\delta} g(y)\frac{\sin [y \mu_n(\pm y; x)]}{y}dy=
$$
$$
    =g(+0)\int\limits_0^{\delta} \frac{\sin [y \mu_n(\pm y; x)]}{y}dy+
    \int\limits_0^{\delta} [g(y)-g(+0)]\frac{\sin [y \mu_n(\pm y; x)]}{y}dy,
$$
then by lemma \ref{L.7} it is sufficient to verify  that the second integral from the right side of this equation tend to zero as $n \to \infty.$

To prove this fact we take $\varepsilon >0,$  choose $h<\delta$ such that
$$
    0\le g(y)-g(+0)<\varepsilon, \quad 0<y\le h,
$$
and divide this integral on two parts:
$$
    i_1 (x;n)+i_2(x;n):=(\int\limits_0^h+\int\limits_h^{\delta})
    [g(y)-g(+0)]\frac{\sin [y \mu_n(\pm y; x)]}{y}dy.
$$

Given the bounded of integrals of kind
$$
    \int\limits_0^h \frac{\sin [y \mu_n(\pm y; x)]}{y}dy, \quad h>0,
$$
which follows from lemma \ref{L.7} and applying the second theorem about average, at each fixed $x \in \mathbb{R}$ we get
\begin{equation}\label{i1}
   | i_1(x;n)|=[g(h)-g(+0)] \bigg| \int\limits_\theta^h \frac{\sin [y \mu_n(\pm y; x)]}{y}dy \bigg| \le \varepsilon C_1(x),
\end{equation}
herewith uniformly over $n \in \mathbb{N}.$

Assume
$$
    \varphi(y):=\cases{0,                        & $0\le y < h,$ \cr
                        \frac{g(y)-g(+0)}{y}, & $h \le y.$}
$$
It is clear that  $\varphi \in L(\mathbb{R}_+).$ Taking this into account, by lemma \ref{L.7'} we  get the estimation
\begin{equation}\label{i2}
    |i_2(x;n)| <\varepsilon C_2(x), \quad n>n_0(\varepsilon).
\end{equation}

From the relations (\ref{i1}) -- (\ref{i2}) implies
$$
    \lim\limits_{n\to \infty} (i_1(x;n)+i_2(x;n))=0,
$$
which proves this lemma. Lemma is proved.

{\bf 6. Proof of the main results.} We will obtain now the main results on  convergence of Fourier series on the system (\ref{system}).

{\bf Proof of Theorem \ref{T.1}.} We show first that for any function $f \in L_p(\mathbb{R}),~ p>1,$ the inequality is true
$$
     \int\limits_{-\infty}^\infty | S_{n}(f;{\bf a};x)|^p dx
     \le  C_p \int\limits_{-\infty}^\infty |f(x)|^p dx.
$$

By lemma \ref{L.5} the representation of the partial sum of the Fourier series on the system $\{\Psi_n(z)\}$ can be written as
\begin{equation}\label{5}
    S_{n}(f;{\bf a};x)=\lim_{\varepsilon\to 0+}
    \frac{1}{\pi} \int\limits_\varepsilon^\infty \frac{g(x+y)-g(x-y)}{y} dy
\end{equation}
where
$$
    g(x\pm y)=f(x\pm y)\sin \big(\pm y \mu_n(\pm y;x)\big),
$$
and the value $\mu_n(\pm y;x)$ is defined by the equality (\ref{def-mu}).

In the book \cite[c. 149]{Tichmarsh_Int_Fourier} was proved that for any function
 $\varphi \in L_p(\mathbb{R}),~ p>1,$ the formula (that called  Hilbert transform)
$$
    H(f;x):=\lim_{\varepsilon \to 0+} \int\limits_{\varepsilon}^\infty
    \frac{\varphi(x+t)-\varphi(x-t)}{t}dt,
$$
almost everywhere determines the function $H(f) \in L_p(\mathbb{R})$ and the inequality holds
$$
     \int\limits_{-\infty}^\infty |H(f;x)|^p dx \le  C_p \int\limits_{-\infty}^\infty |f(x)|^p dx,
$$
where $M_p$ is a positive constant, which depends only on $p.$

Using the fact of bounded of Hilbert transform in the spaces $L_p(\mathbb{R})$ at $p>1,$ given  that the function $g \in L_p(\mathbb{R})$, on the basis of relation (\ref{5}), we get
$$
     \int\limits_{-\infty}^\infty | S_{n}(f;{\bf a};x)|^p dx=
     \int\limits_{-\infty}^\infty | H(g;x )|^p dx\le  C_p \int\limits_{-\infty}^\infty |f(x)|^p dx,
$$
{\it quod erat demonstrandum}.

We now proceed to proof of the assertion of Theorem. Since  $f \in L_p(\mathbb{R}),$
then from results \cite{Kober_1944} it follows that the function $f$ can be  approximated by using fractions (\ref{element-racfunct}) with any degree of accuracy. That is, for any $\varepsilon>0$ can choose a finite linear combination
$$
    T_{n_0}(x)=\sum\limits_{k=0}^{n_0} A_k (x-\bar{a}_k)^{-1} + \sum\limits_{k=1}^{n_0} B_k (x-a_k)^{-1},
$$
such that
$$
    \int\limits_{-\infty}^\infty | f(x)-T_{n_0}(x)|^p dx\le \varepsilon.
$$

 Since the system (\ref{system}) is the result of orthogonalization of the system (\ref{element-racfunct}), then
$$
    T_{n_0}(x)=\sum\limits_{k=0}^{n_0} A_k' \Phi_k(x) + \sum\limits_{k=1}^{n_0} B_k' \Psi_k(x),
$$
and
$$
    S_n(T_{n_0}; x)=T_{n_0}(x), \quad n_0 \le n.
$$

Note this, we get
$$
    \int\limits_{-\infty}^\infty | f(x)-S_{n}(f;x)|^p dx\le
    \int\limits_{-\infty}^\infty | f(x)-T_{n_0}(x)|^p dx +
$$
$$
   + \int\limits_{-\infty}^\infty | S_n(f-T_{n_0};x)|^p dx \le (1+C_p)\varepsilon,
$$
where $\varepsilon$ is arbitrary positive number. Theorem is proved.

{\bf Proof of Theorem \ref{T.2}.} Since function $f$  is summable on $\mathbb{R}$, then from lemma \ref{L.7'} it follows that the representation (\ref{Sn}) can be written in the form
\begin{equation}\label{6}
    S_{n}(f;{\bf a};x_0)=\frac{1}{\pi} \int\limits_0^{\delta} f(x_0-y)
    \frac{\sin y \mu_n(-y;x_0)}{y}~dy+
$$
$$
    +\frac{1}{\pi} \int\limits_0^{\delta} f(x_0+y)
    \frac{\sin y \mu_n(y;x_0)}{y}~dy +o(1), \quad n \to \infty,
\end{equation}
where $\delta>0$ is arbitrary positive number.

By the condition of Theorem  the function $f(\cdot)$ has bounded variation on $\mathbb{R}$, therefore that  function can be represented as a difference of two monotone increasing functions. Applying lemma \ref{L.8} to the each  monotonic component of the function $f$, we get
$$
    \lim\limits_{n\to \infty} S_{n}(f;{\bf a};x_0)= \frac{1}{\pi}\cdot \frac{\pi}{2}[f(x_0-0)+f(x_0+0)]=
    \frac{1}{2}[f(x_0-0)+f(x_0+0)],
$$
{\it quod erat demonstrandum}. Theorem is proved.

{\bf Prof of Theorem \ref{T.3}.}  From Theorem \ref{T.1} implies that
\begin{equation}\label{7}
    \frac{f(x_0-0)+f(x_0+0)}{2}=
    \frac{1}{\pi} \int\limits_0^{\delta} f(x_0-0)
    \frac{\sin y \mu_n(-y;x_0)}{y}~dy+
$$
$$
    +\frac{1}{\pi} \int\limits_0^{\delta} f(x_0+0)
    \frac{\sin y \mu_n(y;x_0)}{y}~dy +o(1), \quad n \to \infty,
\end{equation}
where $\delta>0$ is arbitrary positive number.

From the relations (\ref{6}) and (\ref{7}) we obtain
$$
    S_{n}(f;{\bf a};x_0)-\frac{f(x_0-0)+f(x_0+0)}{2}=
$$
$$
    =\frac{1}{\pi} \int\limits_0^{\delta}
    \frac{f(x_0)-f(x_0-0)}{y}\sin y \mu_n(-y;x_0)~dy+
$$
$$
    +\frac{1}{\pi} \int\limits_0^{\delta}
    \frac{f(x_0+y)f(x_0+0)}{y}\sin y \mu_n(y;x_0)~dy +o(1), \quad n \to \infty.
$$
Since  the  integrals (\ref{T3-cond}) is converged, then statement of the theorem follows from the lemma \ref{L.7'}. Theorem is proved.


\bibliographystyle{plain}
\renewcommand{\refname}{References}

\end{document}